\documentclass[12pt]{amsart}

\usepackage{amsfonts, amssymb, amscd}

\def\dual                 {{\vee}}
\def\rk                 {{\rm rk}}
\def\ii                 {{\rm i}}
\def\ee                 {{\rm e}}

\def\Coker           {{\rm Coker}}
\def\Tor			{{\rm Tor}}

\def\ZZ                 {{\mathbb Z}}
\def\PP                {{\mathbb P}}

\def\CC                 {{\mathbb C}}
\def\QQ                 {{\mathbb Q}}

\newtheorem{lemma}{Lemma}[section]
\newtheorem{theorem}[lemma]{Theorem}
\newtheorem{corollary}[lemma]{Corollary}
\newtheorem{proposition}[lemma]{Proposition}
\theoremstyle{definition}
\newtheorem{definition}[lemma]{Definition}

\newtheorem{remark}[lemma]{Remark}
\theoremstyle{remark}
\newtheorem*{proof*}{Proof}

\title{On the $K$-theory of smooth toric DM  stacks}

\author{Lev A.  Borisov and R. Paul Horja}
\address{Department of Mathematics \\ University of Wisconsin \\
  Madison \\ WI \\ 53706 \\ USA\\{\tt borisov@math.wisc.edu}\\
 The Fields Institute\\ 222 College St\\ 
Toronto\\ Ontario \\ M5T 3J1 \\ Canada
\\{\tt horja@fields.utoronto.ca}.}

\thanks{The first 
author was partially supported by NSF grant DMS-0140172.}

\begin{document}

\begin{abstract}
We explicitly calculate the Grothendieck $K$-theory ring
of a smooth toric Deligne-Mumford stack and define
an analog of the Chern character. In addition, we calculate
$K$-theory pushforwards and pullbacks for weighted blowups
of reduced smooth toric DM stacks.
\end{abstract}

\maketitle

\section{Introduction}
In this paper we calculate the Grothendieck $K$-theory rings
of coherent sheaves on smooth toric Deligne-Mumford stacks.
These stacks, which generalize the notion of a smooth toric
variety, have been defined in \cite{BCS}.

We are mostly interested in the reduced case, which is
characterized by the condition that there is an open
substack which is a subscheme. However, since our technique
is applicable to non-reduced stacks as well, we extend
the result to the general case in a separate section.
We find that the Grothendieck group of a smooth toric
Deligne-Mumford stack $\PP_{\bf\Sigma}$ 
is generated by classes of 
invertible sheaves, and we find the generators of the
ideal of relations satisfied by these sheaves.
In the reduced case $K_0(\PP_{\bf\Sigma})$ is generated
by the classes $R_i$ of the invertible sheaves
${\mathcal L}_i$ which correspond to the one-dimensional
cones of the fan $\Sigma$. 
These sheaves generalize the sheaves ${\mathcal O}(D_i)$
for codimension one strata $D_i$ in a
smooth toric variety. Moreover, we find the generators of the ideal
of relations among $R_i^{\pm 1}$.

\smallskip

\noindent
{\bf Theorem \ref{main}}. 
Let $B$ be the quotient of the Laurent polynomial ring
$\ZZ[x_1,x_1^{-1},\ldots,x_n,x_n^{-1}]$ by the ideal generated 
by the relations
\begin{itemize}
\item
$\prod_{i=1}^n x_i^{f(v_i)}=1$ for any linear function
$f:N\to\ZZ$,
\item
$\prod_{i\in I}(1-x_i)=0$
for any set $I\subseteq [1,\ldots,n]$ such that $v_i,i\in I$
are not contained in any cone of $\Sigma$.
\end{itemize}
Then the map $\rho:B\to K_0(\PP_{\bf\Sigma})$ which sends
$x_i$ to $R_i$ is an isomorphism.

\smallskip

We also show that the $K$-theory with complex coefficients
is a semilocal Artinian $\CC$-algebra and 
explicitly describe its local components.
We produce a \emph{vector space} isomorphism 
between the $K$-theory with complex coefficients
and a combinatorially defined ring dubbed SR-cohomology,
which we call \emph{combinatorial Chern
character}. It generalizes the usual Chern character
in the case of projective toric varieties and is expected to
coincide with the equivariant Chern character 
of \cite{Adem-Ruan} (see also the very recent preprint \cite{JKK})
in the projective DM stack case.

The paper is structured as follows. In Section \ref{sec2}
we recall the definition and basic properties of smooth 
toric DM stacks. We restrict our attention to the reduced 
case, which makes the Gale duality construction of \cite{BCS}
significantly easier to describe. In Section \ref{sec3} 
we introduce the SR-cohomology  in the reduced case and
describe its decomposition into sectors. In Section \ref{sec4}
we calculate the $K$-theory of the reduced smooth toric DM
stack. The key idea is to use the homogenization trick
of \cite{Eisenbud} to resolve any coherent 
sheaf on $\PP_{\bf\Sigma}$
by direct sums of invertible sheaves. Section \ref{sec5} 
defines the combinatorial Chern character, and Section \ref{sec6}
extends the results of Sections \ref{sec4} and \ref{sec5}
to the nonreduced case. Sections \ref{sec7}, \ref{sec8}
and \ref{sec9} describe the $K$-theory pullbacks and pushforwards for special classes of morphisms between
reduced smooth toric DM stacks. These results will be used in
a subsequent paper \cite{BH} on homological mirror symmetry
and the GKZ hypergeometric system of partial differential equations.

There is some overlap between the results of this paper and those of 
the recent preprint \cite{Ba}, where the equivariant 
K-theory of smooth toric varieties is studied with the help of the  
Merkurjev's spectral sequence. 

\section{Review of reduced smooth toric DM stacks}\label{sec2}
In this section we will briefly review the definitions of 
toric Deligne-Mumford stacks, as developed in \cite{BCS}.
We are specifically interested in the reduced case, which
simplifies the construction significantly. 

Let $N$ be a free abelian group, and let $\Sigma$ be 
a simplicial fan in $N$. See \cite{Fulton} for the definition
of a fan. A \emph{stacky fan} ${\bf \Sigma}$ is defined as
the data $(\Sigma,\{v_i\})$ where $\{v_i\}$ is a collection
of lattice points, one for each one-dimensional cone 
$C_i\in\Sigma$. We require $v_i\in C_i\cap (N-\{\bf0\})$, but 
$v_i$ may or may not be the minimum lattice point on $C_i$.
The paper \cite{BCS} makes an additional technical assumption 
\begin{equation}\label{nondegen}
\emph{elements~$v_i$~span~a~finite~index~subgroup~of~$N$}.
\end{equation}

\begin{remark}
It is likely that this assumption \eqref{nondegen}
is just an artifact of the 
technique of \cite{BCS}. In general, one can always 
look at the preimage $N_f$ in $N$ of the torsion part
of $N/Span(v_i,i=1,\ldots,n)$. The quotient $N/N_f$ splits
off (noncanonically) as a free direct summand. 
It is reasonable to expect that the DM stack that 
corresponds to $\bf\Sigma$ is the product of 
$(\CC^*)^{{\rm rk}(N/N_f)}$ and the toric DM stack of \cite{BCS}
constructed for ${\bf \Sigma}$ in $N_f$ instead of $N$.
However, the  functoriality of this construction is a bit unclear,
in view of the fact that the splitting is non--canonical.
\end{remark}

The collection $\{v_i\},i=1,\ldots,n$ gives a map 
$$
\phi:\ZZ^n\to N
$$
with finite cokernel. We dualize to get
an injective map 
$$
\phi^\dual:N^\dual\to (\ZZ^n)^\dual.
$$
The Gale dual of $\phi$ (see \cite{BCS}) is simply the map 
$
(\ZZ^n)^\dual \to \Coker(\phi^\dual).
$
Denote by $G$ the algebraic group 
$
{\rm Hom}(\Coker(\phi)^\dual,\CC^*).
$
The exact sequence
\begin{equation}\label{cokersurj}
0\to N^\dual\to(\ZZ^n)^\dual\to \Coker(\phi^\dual)\to 0
\end{equation}
leads to the exact sequence
$$
0\to G \to(\CC^*)^n\to {\rm Hom}(N^\dual,\CC^*).
$$
Hence, $G$ can be thought of as a subgroup of $(\CC^*)^n$.
Its $\CC$-points can be 
thought of as collections of $n$ nonzero complex numbers
$\lambda_i, i=1,\ldots,n$ which satisfy the condition
\begin{equation}\label{eqlam}
\prod_{i=1}^n \lambda_i^{f(v_i)}=1
\end{equation}
for all linear functions $f:N\to \ZZ$.

Consider the subset $Z$ of $\CC^n$ that consists of all the
points ${\bf z}=(z_1,\ldots,z_n)$ such that the set of $v_i$ for
the zero coordinates of ${\bf z}$ is contained in a cone of $\Sigma$.
Then the toric DM stack $\PP_{\bf\Sigma}$
that corresponds to the stacky fan
${\bf\Sigma}=(\Sigma,\{v_i\})$ is defined as the stack quotient 
$[Z/G]$ where $Z$ and $G$ are 
endowed with the natural reduced 
scheme structures. The action of $(\lambda_i)$ on
$(z_i)$ is given by $(\lambda_iz_i)$.
It has been shown in \cite{BCS} that $\PP_{\bf\Sigma}$
is a Deligne-Mumford stack whose moduli space is
the simplicial toric variety $\PP_\Sigma$.

To a cone in $C$ one can associate a closed substack of 
$\PP_{\bf\Sigma}$ by looking at the quotient of $N$ by the sublattice
spanned by $v_i\in C$. The new fan is defined as the image of the 
link of $C$ in $\Sigma$, and the new $v_j$ are the images
of the old $v_j$ from the link. Unfortunately, this procedure
may introduce torsion into $N$, so the resulting substack
is not reduced. Moreover, in order for the construction of
\cite{BCS} to work, the images of $v_j$ from the link of 
$C$ need to generate a finite index subgroup in 
$N/Span(v_i\in C)$. The easiest way to assure this is 
by imposing a stronger technical condition
\begin{equation}\label{technical}
\emph{every~cone~of~$\Sigma$~is~contained~in~a~
cone~of~dimension~{\rm rk}$N$},
\end{equation}
see \cite{BCS}. This also assures that $\PP_{\bf\Sigma}$ is 
covered by open substacks of the form $\CC^{{\rm rk}N}/G_i$
where $G_i$ are finite abelian groups. These open substacks
can be indexed by the cones in $\Sigma$ of maximum dimension.
On the other hand, one can definitely talk about closed
toric subvarieties of the coarse moduli space $\PP_\Sigma$
of $\PP_{\bf\Sigma}$ without the additional assumption 
\eqref{technical}. Throughout this paper we will only use
the assumption \eqref{nondegen}.

\section{Orbifold cohomology and SR-cohomology of reduced smooth toric DM stacks}\label{sec3}
\label{secSR}
In this section  we introduce
some combinatorial invariants of toric DM stacks which 
we call SR-cohomology rings. They coincide
with the orbifold cohomology rings in the projective case
but are generally different, even for smooth toric varieties.

A natural combinatorial invariant of a fan is the partial 
semigroup ring $\ZZ[N,\Sigma]$ which is defined 
as a free abelian group with the
basis $\{[w],w\in N\cap\Sigma\}$ indexed by 
the lattice elements inside the support of the fan and multiplication
$$
{[w_1]}*{ [w_2]}=\Big\{
\begin{array}{ll}
[{w}_1+{w}_2],&{\rm if~there~exists~} C\in\Sigma,
C\ni{w}_1,{ w}_2\\
0,&{\rm otherwise}.
\end{array}
$$
The rings $\QQ[N,\Sigma]$ and $\CC[N,\Sigma]$ are 
defined analogously. For a stacky fan $\bf\Sigma$,
these rings are given additional structure of graded rings.
The grading can take nonnegative rational values.  It
is defined by $\deg([w])=\sum_i{\alpha_i}$ where 
$w=\sum_{i}\alpha_i v_i $ is the unique way of writing
$w$ as a rational linear combination of $v_i$ that lie 
in the minimum cone of $\Sigma$ that contains it.

To any stacky fan $\bf\Sigma$ one can associate a
\emph{SR-cohomology ring}, to be denoted by $H_{SR}(\PP_{\bf\Sigma},\CC)$.
\begin{definition}
Pick a basis $\{f_1,\ldots,f_{\rk N}\}$ of ${\rm Hom}(N,\ZZ)$.
Each $f_i$ gives rise to a degree one element 
$t_i:=\sum_{j=1}^n f_i(v_j)[v_j]$ in 
$\CC[N,{\Sigma}]$.
Then 
$$
H_{SR}(\PP_{\bf\Sigma},\CC) :=  
\CC[N,\Sigma]/\langle t_1,\ldots,t_{\rk N}\rangle.
$$
The rings $H_{SR}(\PP_{\bf\Sigma},\QQ)$ and $H_{SR}(\PP_{\bf\Sigma},\ZZ)$
are defined analogously.
\end{definition}

It is clear that the above defined SR-cohomology
rings are independent of the choice of the basis $\{f_i\}$. If $\PP_{\bf\Sigma}$
is a smooth projective toric \emph{variety}, then 
its combinatorial cohomology rings are isomorphic to
its usual cohomology rings. Indeed, this is precisely the 
Stanley-Reisner presentation of the cohomology ring of 
a smooth projective toric variety. This is the motivation
behind our notation $H_{SR}$. We 
call our rings SR-cohomology rings, as opposed
to Stanley-Reisner cohomology rings, to avoid confusion
with the term
Stanley-Reisner ring that is already used in the literature.

\begin{remark}
We abuse notation slightly and still use 
$H_{SR}(\PP_{\bf\Sigma},*)$ even if $\bf\Sigma$
does not satisfy \eqref{nondegen}. Hopefully, one
will eventually define toric stacks without this restriction.
\end{remark}

The main technical result of \cite{BCS}
connects the SR-cohomology of a smooth toric
DM stack with its orbifold cohomology, see \cite{CR} 
and 
\cite{AGV}. However, one needs an important 
additional assumption 
that the coarse moduli space $\PP_\Sigma$ is projective.
\begin{theorem}(\cite{BCS})
For any $\PP_{\bf\Sigma}$  such that $\PP_\Sigma$ is
projective, there holds
$$
H_{orb}(\PP_{\bf\Sigma},\QQ)\cong H_{SR}(\PP_{\bf\Sigma},\QQ).
$$
\end{theorem}

\begin{remark}\label{notsofast}
In general, the SR-cohomology and orbifold or usual 
cohomology rings are
different. The simplest example is given by the fan $\Sigma$
which is the union of the first and the third quadrants in
$\ZZ^2$. The elements $v_i$ are $(\pm 1,0)$ and $(0,\pm 1)$.
The corresponding variety is $\PP^1\times\PP^1
-\{(0,\infty),(\infty,0)\}$. The SR-cohomology is three-dimensional,
with dimension $2$ graded component of degree one.
It is unclear what geometrically defined cohomology theory 
on $\PP^1\times\PP^1
-\{(0,\infty),(\infty,0)\}$ can produce such a ring.
\end{remark}

Despite Remark \ref{notsofast}, the SR-cohomology ring $H_{SR}(\PP_{\bf\Sigma},\CC)$ is very 
well behaved. 
\begin{proposition}
The SR-cohomology ring $H_{SR}(\PP_{\bf\Sigma},\CC)$ 
of  any (reduced) stacky fan ${\bf \Sigma}$
is a finite dimensional complex vector space. It is a
local Artinian graded $\CC$-algebra with the maximum ideal 
given by the span of the elements of positive degree.
\end{proposition}

\begin{proof}
Since $H_{SR}(\PP_{\bf\Sigma},\CC)$
is nonnegatively graded and its degree zero part 
is isomorphic to $\CC$, it is enough to show that it is
an Artinian ring. Hence, let us study 
$\CC$-algebra homomorphisms 
$\phi:H_{SR}(\PP_{\bf\Sigma},\CC)\to \CC$.

In view of relations in $\CC[N,\Sigma]$, the
values of $\phi([w])$ are nonzero only for $w$ in
some cone $C$ of $\Sigma$. Then the relations
$t_i$ show that these values must be zero
for $w\in C$ that are integer linear combinations of
$v_j\in C$. Finally, for any other nonzero $w\in C$, some
positive multiple of it is an integer linear combination
of $v_j\in C$. This leads to $\phi([w])^l=0$, so
$\phi([w])=0$ for all nonzero $n$. This shows that
the only maximum ideal in $H_{SR}(\PP_{\bf\Sigma},\CC)$ 
is the span of the 
elements of positive degree.
\end{proof}

In the case of projective $\PP_\Sigma$, there is 
a natural decomposition of the orbifold and 
SR-cohomology into direct sum of sectors.
Each sector corresponds to the usual cohomology of some
closed toric subvariety of $\PP_\Sigma$.
We observe that this decomposition still occurs for the 
SR-cohomology, even without the projectivity assumption. 
Moreover, we will not even assume \eqref{nondegen}.
\begin{definition}
The \emph{untwisted sector} of SR-cohomology
is defined as the subring of $H_{SR}(\PP_{\bf\Sigma},\CC)$
 generated by the images of $[v_i],i=1,\ldots,n$.
\end{definition}

\begin{remark}
It is easy to see that the choice of $v_i$ in the 
corresponding one-dimesnional cone 
does not 
change the untwisted sector much, if one works with
rational or complex coefficients. Indeed, one can
simply rescale the corresponding variables. Consequently,
we can talk about the SR-cohomology
of the toric variety $\PP_\Sigma$, in analogy with the 
projective case. We can define it as the untwisted
sector of the SR-cohomology 
of the stack $\PP_{\bf\Sigma}$ obtained by picking 
the $v_i$ to be the smallest lattice points on the one-dimensional
cones of $\Sigma$. We  denote these rings by
$H_{SR}(\PP_\Sigma,\QQ{~\rm or~}\CC)$.
\end{remark}

For a cone $C\in\Sigma$, let
${\rm Box}(C)$ be
the set of elements $v$ of $N$ which are linear
combinations with rational coefficients in the range $[0,1)$
of elements $v_i\in C$.  
Let ${\rm Box}({\bf\Sigma})$ be
the union of ${\rm Box}(C)$ for all $C\in\Sigma$.

\begin{proposition}\label{split}
As a module over the untwisted 
sector, the ring
$H_{SR}(\PP_{\bf\Sigma},\CC)$ is a direct sum 
of the modules $H_{v}$ generated by the images 
in $H_{SR}(\PP_{\bf\Sigma},\CC)$ of the elements 
$[v]\in\CC[N,\Sigma]$ for all $v\in{\rm Box}({\bf\Sigma})$. 
\end{proposition}

\begin{proof}
The statement follows from the analogous result 
for $\CC[N,\Sigma]$, where it is obvious.
\end{proof}

In the projective case, each of the modules $H_v$
is isomorphic to the cohomology of the  toric subvariety 
of $\PP_\Sigma$ that corresponds to the minimum
cone of $\Sigma$ that contains $v$. Remarkably, this 
is true in general
if one works in SR-cohomology with complex or rational
coefficients.
\begin{proposition}\label{sector}
The module $H_v$ is isomorphic to the SR-cohomology 
with complex coefficients of the (closed) toric subvariety 
of $\PP_\Sigma$ that corresponds to the minimum
cone of $\Sigma$ that contains $v$.
\end{proposition}

\begin{proof}
The proof is analogous to that of \cite[Proposition 5.2]{BCS},
but we briefly sketch it here for the benefit of the reader.
Let $N_1$ denote the quotient of $N$ by the subgroup generated
by $v_i\in C$, and let ${N'}=N_1/torsion$ be
its torsion-free part. Let $\Sigma'$ in $N'$ be the image
of the link of $C$ in $\Sigma$ and let $v_i'$ be the images
of $v_i$ from the link.
The SR-cohomology of the closed
toric subvariety that corresponds to $C$ is isomorphic 
to the untwisted sector in the SR-cohomology 
$H_{SR}(\PP_{{\bf \Sigma}'},\CC)$ of ${\bf \Sigma}'$. 

Observe that $H_v$ is a quotient of the submodule $M_v$
of 
$\CC[N,\Sigma]$ which is supported on the star of $C$.
Moreover, it is generated by the monomial of the form
$[n+v]$ where $n$ is a lattice point in the star of $C$ in
$\Sigma$ which is an \emph{integer} linear combination 
of the $v_j$  in the minimum cone that contains it.
This identifies $H_v$ with the quotient of the polynomial ring 
in the variables $D_1,\ldots, D_k$ that correspond to 
one-dimensional faces of $C$ and cones in its link
by the ideal with generators
\begin{itemize}
\item
$\prod_{i\in I}D_i$, if no cone $C'\supseteq C$
contains all $v_i,i\in I$,
\item
$\sum_{i,v_i\in {\rm star}(C)} f(v_i)D_i$
for any linear map $f:N\to \ZZ$.
\end{itemize}

On the other hand, the untwisted sector of 
$H_{SR}(\PP_{{\bf \Sigma}'},\CC)$ is
isomorphic to the quotient of the polynomial ring in 
the variables $D'_1,\ldots, D'_{k'}$ which correspond to
one-dimensional faces of cones in the link of $C$ 
(hence not in $C$), by the ideal generated by 
\begin{itemize}
\item
$\prod_{i\in I}D'_i$, if no cone $C'\supseteq C$
contains all $v_i,i\in I$,
\item
$\sum_{i,v_i\in {\rm link}(C)} f'(\rho(v_i))D'_i$
for any linear map $f':N'\to \ZZ$.
\end{itemize}
Here $\rho:N\to N'$ is the projection.

To connect these two spaces, observe that 
linear functions on $N'$ lift to linear functions on $N$.
Pick a complementary sublattice in the lattice of linear
functions on $N$ and pick its basis $f_1,\ldots,f_{\dim C}$.
These $f_i$ provide relations on $D_i$ that 
allow one to express $D_i,v_i\in C,$ in terms of 
$D_i, v_i\in{\rm link}(C),$ in the relations for $H_v$.
The remaining relations are precisely those 
of the untwisted sector of $H_{SR}
(\PP_{{\bf\Sigma}'},\CC)$. Moreover, this isomorphism is
independent of the choice of the complementary lattice
or its basis.
\end{proof}

\begin{remark}
Propositions \ref{split} and \ref{sector} still hold for SR-cohomology with 
rational coefficients. The proofs are unchanged.
\end{remark}

\begin{remark} Suppose that the fan $\Sigma$ is projective, or 
that it is a subdivision of a cone. Under these assumptions,
if one picks a basis $\{f_j\}$ of $N^\dual$, then the corresponding
elements $\sum_{i=1}^n f(v_i)[v_i]$ form a regular 
sequence in $\QQ[N]^{\Sigma}$. As a result, the 
graded dimension of $H_{SR}({\bf\Sigma},\CC)$, defined
as $\sum_{d\in\QQ}t^d\dim_\CC
H_{SR}(\PP_{\bf\Sigma},\CC)_{\deg=d}
$, equals
$$ 
{\rm gr.dim}H_{SR}(\PP_{\bf\Sigma},\CC)=
(1-t)^{\rk N}\sum_{{\bf n}\in N\cup\Sigma}t^{\deg({\bf n})}. 
$$
The proof of this regularity is 
sketched in \cite{BCS} for the projective case and in \cite{BM}
for the cone case. In general, the graded dimension of
$H_{SR}({\bf\Sigma},\CC)$ is the sum over sectors 
of the graded dimensions of the SR-cohomology of the sector,
but it lacks such a nice combinatorial formula.
\end{remark}

\section{$K$-theory of reduced toric DM stacks}\label{sec4}
The goal of this section is to prove a combinatorial description
 of the
$K$-theory of reduced toric DM stack which is analogous
to the Stanley-Reisner presentation of the cohomology
of a smooth toric variety. Analogous statements for
smooth toric varieties are contained in \cite{Ba}.
The resulting ring is then compared to the SR-cohomology
of the stack. 
We use the notations from the previous sections.
We do not make any assumptions on the stacky fan
apart from
\eqref{nondegen}.

\begin{definition}
Let $\mathcal X$ be a smooth Deligne-Mumford stack.
Define the (Grothendieck) $K$-theory group 
$K_0({\mathcal X})$ to be the quotient of the 
free abelian group generated by coherent 
sheaves $\mathcal F$ on $\mathcal X$
by the relations 
$[{\mathcal F}_1]-[{\mathcal F}_2]+[{\mathcal F}_3]$ 
for all exact sequences
$0\to{\mathcal F}_1\to {\mathcal F}_2\to {\mathcal F}_3\to 0$.
\end{definition}

\begin{remark}
The $K$-theory of $\mathcal X$ admits a product structure by
$$
[{\mathcal F}_1]*[{\mathcal F_2}]=\sum_{i=0}^{\dim {\mathcal X}}
(-1)^i[{\mathcal T}\hskip -.1cm
or^i({\mathcal F}_1,{\mathcal F}_2)].
$$
The image of the 
structure sheaf $\mathcal O$ plays the role of the 
identity.
\end{remark}



We recall that the category of coherent sheaves on $[Z/G]$ 
is equivalent to that of $G$-linearized coherent sheaves on
$Z$, see \cite[Example 7.21]{Vistoli}. We will always implicitly
use this equivalence.
\begin{definition}
For each $i$ from $1$ to $n$ we define a 
$G$-linearized invertible sheaf
${\mathcal L}_i$ on $Z$ as follows.
As a sheaf, it will be isomorphic to the structure sheaf
${\mathcal O}_Z$. For $g=(\lambda_1,\ldots,\lambda_n)
\in G$ the isomorphism 
${\mathcal O}_Z\to g^*{\mathcal O}_Z={\mathcal O}_Z$ sends
$1$ to $\lambda_i$.
Here we have used the canonical isomorphism $g^*{\mathcal O}_Z={\mathcal O}_Z$.
We define $R_i$ to be the image of ${\mathcal L}_i$
in the $K$-theory  of $\PP_{\bf\Sigma}$.
\end{definition}

\begin{remark}\label{sections}
The sheaf ${\mathcal L}_i$ 
has a $G$-invariant global section given by 
the restriction to $Z$ of the $i$-th coordinate function $z_i$
on $\CC^n$.  Indeed, $g^*z_i = \lambda_i z_i$,
so $z_i$ is compatible with the isomorphisms of the above
definition.
In the case of a smooth toric \emph{variety}, 
${\mathcal L}_i$ corresponds to 
${\mathcal O}(D_i)$ where $D_i$ is
the divisor of the  $i$-th dimension one 
cone of $\Sigma$. 
\end{remark}

\begin{proposition}
The 
elements $R_i$ of $K_0(\PP_{\bf\Sigma})$ satisfy the 
following relations.
\begin{itemize}
\item
$\prod_{i=1}^n R_i^{f(v_i)}=1$ for any linear function
$f:N\to\ZZ$,
\item
$\prod_{i\in I}(1-R_i)=0$
for any set $I\subseteq [1,\ldots,n]$ such that the $v_i,i\in I$
are not contained in any cone of $\Sigma$.
\end{itemize}
\end{proposition}

\begin{proof}
The 
relations  \eqref{eqlam} and the definition of 
${\mathcal L}_i$
show that $\otimes_i {\mathcal L}_i^{f(v_i)}$ 
are in fact trivially linearized on $Z$ for every
linear function $f:N\to \ZZ$. Passing to $K$-theory,
we get the first set of relations on $R_i$.

To get the second set of relations, consider the 
Koszul complex of the set of elements 
$z_i,i\in I$ in the ring $A=\CC[z_1,\ldots,z_n]$.
It gives an exact sequence  
$$
 \ldots \to \oplus_{J\subseteq I,\vert J\vert = k}A\to
 \oplus_{J\subseteq I,\vert J\vert = k-1}A\to \ldots
 $$$$\ldots
\to \oplus_{j\in I} A \to A \to A/\langle z_i,i\in I \rangle\to 0
$$
of $A$-modules. We can pass to the associated sheaves 
on $\CC^n$ and then restrict them to $Z$ to get 
$$
0\to\ldots \to \oplus_{J\subseteq I,\vert J\vert = k}
{\mathcal O}_{Z}
\to
\oplus_{J\subseteq I,\vert J\vert = k-1}
{\mathcal O}_{\CC^n}
\to\ldots\to 0
$$
where we have used the condition on $I$ to see that 
$A/\langle z_i,i\in I \rangle$ is supported outside of $Z$.
The maps between the copies of ${\mathcal O}_Z$
for $J_1$ and $J_2$ are zero unless $J_2\subset J_1$
and are $\pm z_{J_1-J_2}$ otherwise. We can make this 
into a complex of $G$-linearized sheaves 
$$
0\ldots \to \oplus_{J\subseteq I,\vert J\vert = k}
\otimes_{j\in J} {\mathcal L}^{-1}_j
\to
\oplus_{J\subseteq I,\vert J\vert = k-1}
\otimes_{j\in J} {\mathcal L}^{-1}_j
\to\ldots\to 0
$$
by  using the $G$-equivariant maps $z_j: {\mathcal L}^{-1}_j\to
{\mathcal O}_Z$ which are twists by ${\mathcal L}^{-1}_j$
of the maps of Remark \ref{sections}.

The alternating sum of any long exact 
sequence of sheaves is zero in the $K$-theory, which implies
the second set of relations, after multiplying by 
the invertible element $\prod_{i\in I} R_i$.
\end{proof}

\begin{theorem}\label{surj}
The 
elements $R_i$ generate the $K_0$-theory of the reduced toric 
DM stack $\PP_{\bf\Sigma}$.
\end{theorem}

\begin{proof}
Consider a $G$-linearized coherent sheaf $\mathcal F$ on
the open subset $Z$ of $\CC^n$. Let us denote the embedding
of $Z$ into $\CC^n$ by $i$. We observe that 
$i^*i_* {\mathcal F}\to\mathcal F$ 
is an isomorphism, simply because $i$ is an open embedding.
Since $\CC^n$ is affine, $i_*{\mathcal F}$ is the sheaf
associated to $\CC[z_1,\ldots,z_n]$-module 
$M=H^0(Z,{\mathcal F})$. 

\begin{lemma}\label{fingen}
In the notations above, $M$ is a finitely generated module. 
\end{lemma}

\emph{Proof of Lemma \ref {fingen}.}
We know by \cite[Exercise II.5.15]{Hartshorne} that 
$\mathcal F$ is a restriction to $Z$ of some coherent sheaf
on $\CC^n$. This coherent sheaf can be resolved by free
sheaves on $\CC^n$.  When we restrict to $Z$ we see that
$\mathcal F$ is resolved by direct sums of ${\mathcal O}_Z$. 
Hence, it is enough to see that 
all the cohomology  groups of ${\mathcal O}_Z$ 
are finitely generated modules 
over $\CC[z_1,\ldots,z_n]$. By using  the $(\CC^*)^n$-action,
we see that $H^0(Z,{\mathcal O}_Z)$ is generated by 
some monomials $\prod_j z_j^{r_j}$, if one thinks of it as 
the subspace of the quotient field of $\CC[z_1,\ldots,z_n]$.
We easily see that $\CC^n-Z$ is a union of coordinate 
subspaces of codimension at least two, which implies that 
$r_i$ have to be nonnegative. Then we see that 
$H^0(Z,{\mathcal O}_Z)=\CC[z_1,\ldots,z_n]$.
To calculate the higher cohomology $H^*(Z,{\mathcal O}_Z)$ 
groups,
we can cover $Z$ by open affine subsets $U_C$ of the form
$$
U_C = \{(z_1,\ldots,z_n), z_i\neq 0 ~{\rm for}~v_i\not\in C\}
$$
where $C$
runs over the set of nonzero
cones of $\Sigma$. The sections of ${\mathcal O}_Z$
on $U_C$ are the linear span of monomials
$\prod_j z_j^{r_j}$ such that the set of $i$ with $r_i<0$
is contained in the set of $i$ with $v_i\not\in C$. Various 
intersections $U_{C_1C_2\ldots C_k}$
of $U_C$ give the spans of monomials with the condition
that for each $j$ with $r_j<0$ the lattice element $v_j$ is 
not contained in $\cap_lC_l$. We use the cover by sets 
of type $U_C$
to calculate the cohomology of ${\mathcal O}_Z$ as
\v{C}ech cohomology.
For a given $(r_1,\ldots,r_n)$,
the cohomology of ${\mathcal O}_Z$ that corresponds to this
$\ZZ^n$-grading is given by the reduced homology of the following
simplicial complex. Its set of vertices is the set $\Sigma_+$ 
of all cones of positive dimension in $\Sigma$, and its maximum
simplices are the complements in $\Sigma_+$ 
of the one-element subsets $\{C\}$ that correspond to cones
that do not contain any $v_i$ with $r_i<0$. So we have 
a simplicial complex which is a union of a collection of
complements of one-dimensional subsets.  If not all the elements of 
$\Sigma_+$ are in this collection, then the resulting complex
is a cone, since adding such an element has no effect on 
whether a subset of $\Sigma_+$ belongs to the complex.
Hence, it has trivial reduced homology. If all elements
of $\Sigma_+$ are a part of the collection, then 
there can be no $v_i$ with $r_i<0$. The resulting 
complex is a sphere, and we have a one-dimensional top
reduced homology that corresponds to the monomial
in $H^0(Z,{\mathcal O}_Z)$. This shows that the higher cohomology
$H^{>0}(Z,{\mathcal O}_Z)$ groups vanish, which finishes the proof
of the lemma. $\hfill\Box$

\emph{Proof of Theorem \ref{surj} continues.}
The $G$-linearization on $\mathcal F$ gives rise to a $G$-action
on $M$, which is compatible with the $G$-action on
$A=\CC[z_1,\ldots,z_n]$, in the sense that $g(rm)=g(r)g(m)$
for all $g\in G$, $r\in A$ and $m\in M$. Moreover,
we claim that $M$ is generated by a finite set of \emph{eigenelements} 
$m_j$ with $g(m_j)=\chi_j(g)m_j$ for some character 
$\chi_j:G\to\CC^*$. In view of the structure of $G$,
all of its algebraic finite-dimensional actions are diagonalizable.
As a consequence, a module $M$ is generated by a finite set
of eigenelements if and only if 
\begin{equation}\label{fincond}
\begin{array}{c}\emph{for any $m\in M$ the linear 
span of $gm,g\in G$}\\ \emph{is finite-dimensional,
and the action of
$G$ on it is algebraic.}  
\end{array}
\end{equation}
It is clear that this finiteness condition \eqref{fincond}
for a $G$-equivariant Noetherian $A$-module implies
\eqref{fincond} for all equivariant submodules and quotients.

Recall (see \cite{GIT}) that the $G$-linearization
of $\mathcal F$ on $Z$ is given by an isomorphism
$\phi:\mu^* {\mathcal F}\to \pi_2^* {\mathcal F}$ 
on $G\times Z$ where $\mu$ is
the multiplication and $\pi_2$ is the second projection.
Denote by $A_G$ the ring of regular functions on $G$.
We will use the fact that 
$H^0(Z,{\mathcal O}_Z)=H^0(\CC^n,{\mathcal O}_{\CC^n})=A$
from the proof of Lemma \ref{fingen}.
The isomorphism $\phi$ induces a map
$$
\phi_{global}:M\otimes_{A}(A_G\otimes_\CC A)
\to M\otimes_\CC A_G
$$
where the tensor multiplication on the left is via the map
$A\to A_G\times A$ induced by the action of $G$ on $Z$.
Indeed, the left hand side maps to the global sections 
of $\mu^*({\mathcal F})$, whereas the right hand side
consists of the global sections of $\pi_2^* {\mathcal F}$.
This map $\phi_{global}$ encodes the action of $G$ on $M$
by mapping $m\otimes (f\otimes a) 
\mapsto \sum_{k}m_k\otimes f_k$
such that for any $g\in G$ there holds
$a\,g(m)\,f(g)=\sum_k f_k(g) m_k$.
By taking $a=f=1$ we see that
there is a finite sum $\sum_k  m_k\otimes f_k$
such that $g(m)=\sum_k f_k(g) m_k$ for all $g\in G$.
In particular, the span of $gm,g\in G$ is finite-dimensional.
By picking a basis of it, and applying the above, we see
that the action of $G$ is algebraic, which shows
\eqref{fincond}.

Consequently, there is a presentation
$$
F_1\to F_0 \to M\to 0
$$
where $F_i$ is a direct sum of rank one 
$\CC[z_1,\ldots,z_n]$-modules generated
by eigenelements of $G$. Indeed, the kernel of $F_0\to M$
also satisfies \eqref{fincond}. 
Hence, it is generated by a finite number
of eigenelements, which allows one to construct $F_1$.
We will now use the 
homogenization trick of \cite[Corollary 19.8]{Eisenbud}.
Namely, consider the ring $\CC[z_0,\ldots,z_n]$
and extend the $G$-action on it by $gz_0=z_0$ for all $g\in G$.
The map $F_1\to F_0$ is given by a matrix $S$ of polynomials
in $z_1,\ldots,z_n$. Define the matrix $\tilde S$ of homogeneous
polynomials of some large fixed degree $d$ by multiplying each
monomial in $S$ by an appropriate power of $z_0$.
We will then have a presentation of some module $\tilde M$
over $\CC[z_0,\ldots,z_n]$
$$
\tilde F_1\to \tilde F_0\to \tilde M\to 0.
$$
Here $\tilde F_i$ are direct sums of rank one $G$-equivariant
$\CC[z_0,\ldots,z_n]$-modules generated by 
homogeneous eigenelements.
As in \cite[Corollary 19.8]{Eisenbud}, we observe that 
$$
\CC[z_1,\ldots,z_n]=\CC[z_0,\ldots,z_n]/(1-z_0)
$$
and $\tilde M/(1-z_0)\tilde M = M$. Given a $G$-equivariant
homogeneous module $\tilde F$ over $\CC[z_0,\ldots,z_n]$,
we can look at the finite-dimensional vector space 
$$
\tilde F/<z_0,\ldots,z_n>\tilde F.
$$
This vector space inherits the $G$-action. It also inherits
the  the grading,
which is preserved by the $G$-action. As a result, this vector
space has a basis
of homogeneous eigenelements. 
Each such eigenelement can be lifted to
an element $r_{j}$ of $\tilde F$. By looking at the 
(finite-dimensional) $G$-span of $r_{j}$ in $\tilde F$, we can modify $r_{j}$ to be itself a homogeneous eigenelement. 
The elements $r_{j}$ generate $\tilde F$ 
by the graded Nakayama lemma.

We apply this procedure to the kernel $K$ 
of $\tilde F_1\to \tilde F_0$,
and continue on to build a resolution of $\tilde M$. Since
the end of this resolution is the minimal graded resolution
of the homogenous module $K$, it terminates.
We then get a free resolution 
$$
0\to \tilde F_l\to\ldots\to\tilde F_0\to \tilde M\to 0
$$
of $\tilde M$ such that each $\tilde F_i$ is freely generated by
homogeneous eigenelements $r_{ij}$ and the maps are
$G$-equivariant and are compatible with the grading. 
We mod out this resolution by 
$1-z_0$ to get a resolution of $M$. The exactness follows
from $\Tor^{>0}_{\CC[z_0,\ldots,z_n]}(M,\CC[z_1,\ldots,z_n])=0$,
as in \cite[Corollary 19.8]{Eisenbud}. We thus get a 
$G$-equivariant resolution of $M$ by direct sums of 
free modules of rank one generated by eigenelements.
We claim that each such module can be identified with the  
global sections of a tensor product of powers of line bundles
${\mathcal L}_i$. Indeed, one simply needs to show that
every character of $G$ is a restriction of a character of
$(\CC^*)^n$, which follows from the exact sequence \eqref{cokersurj}.

We pass to the corresponding exact sequence
of $G$-linearized sheaves on $\CC^n$ and then restrict it to $Z$
to see that $\mathcal F$ is a linear combination of products
of $R_i$. 
\end{proof}

\begin{corollary}\label{cor:db}
The bounded derived category of the category of coherent
sheaves
on $\PP_{\bf\Sigma}$ is generated by the invertible sheaves
$\otimes_{i}{\mathcal L}_i^{r_i}$.
\end{corollary}

\begin{proof}
The argument of Theorem \ref{surj} shows that every
sheaf on $\PP_{\bf\Sigma}$ admits a free resolution
by direct sums of invertible sheaves of the above type.
\end{proof}

\begin{remark}
We refer the reader to \cite{Ka} for much stronger results
concerning these derived categories.
\end{remark}

\begin{theorem}\label{main}
Let $B$ be the quotient of the Laurent polynomial ring
$\ZZ[x_1,x_1^{-1},\ldots,x_n,x_n^{-1}]$ by the ideal generated 
by the relations
\begin{itemize}
\item
$\prod_{i=1}^n x_i^{f(v_i)}=1$ for any linear function
$f:N\to\ZZ$,
\item
$\prod_{i\in I}(1-x_i)=0$
for any set $I\subseteq [1,\ldots,n]$ such that $v_i,i\in I$
are not contained in any cone of $\Sigma$.
\end{itemize}
Then the map $\rho:B\to K_0(\PP_{\bf\Sigma})$ which sends
$x_i$ to $R_i$ is an isomorphism.
\end{theorem}

\begin{proof}
Theorem \ref{surj} shows that $\rho$ 
is surjective. It is therefore sufficient to show its injectivity.
We will define a map $\rho_1:K_0(\PP_{\bf\Sigma})\to B$
and prove that $\rho_1\circ\rho={\bf id}$.

For any $G$-linearized sheaf $\mathcal F$ on $Z$,
consider the $G$-equivariant module 
$M=H^0(Z,{\mathcal F})$ over $A=\CC[z_1,\ldots,z_n]$.
Consider the $A$-module $\CC$ with trivial $G$-action.
For each $i$ from $0$ to $n$ the finite-dimensional
vector space $\Tor_A^i(M,\CC)$ is acted upon by $G$.
It is a direct sum over the characters of $\chi:G\to \CC^*$ 
of the eigenspaces $V_{i,\chi}$. The group ${\rm Hom}(G,\CC^*)$
is a quotient of $\ZZ^n$ by $N^\dual$. 
Because of the first set of relations on $x_i$, every character
$\chi$ gives a well-defined monomial $x_\chi \in B$. 
We then define 
$$
\rho_1:\mathcal F\mapsto \sum_\chi\sum_{i=0}^n (-1)^i \dim_\CC
(V_{i,\chi}) x_\chi.
$$

We would like to show that $\rho_1$ passes down to
a well-defined map on $K$-theory, i.e. it is additive on short
exact sequences. If 
\begin{equation}\label{exact}
0\to{\mathcal F}_1\to{\mathcal F}_2\to {\mathcal F}_3\to 0
\end{equation}
is an exact sequence of $G$-linearized sheaves on $Z$
with $G$-equivariant morphisms, then the sequence
of $A$-modules 
$$
0\to M_1\to M_2\to M_3
$$
is only exact on the left. We complete it to 
a long exact sequence
$$
0\to M_1\to M_2\to M_3 \to M_4 \to 0.
$$
This long exact sequence of modules splits into short 
exact sequences, which, in turn, give long exact 
sequences of $\Tor^i(*,\CC)$. This shows that 
$$\rho_1(M_1)+\rho_1(M_3)-\rho_1(M_2)=
\sum_\chi\sum_{i=0}^n (-1)^i\dim_\CC(\Tor^i(M_4,\CC)_\chi)
R_\chi.
$$

Because of \eqref{exact}, $M_4$ is associated to the 
$G$-equivariant sheaf ${\mathcal F}_4$ which is 
supported on $\CC^n-Z$. 
We will show that 
\begin{equation}\label{zero}
\sum_\chi\sum_{i=0}^n (-1)^i\dim_\CC(\Tor^i(M_4,\CC)_\chi)
x_\chi=0
\end{equation}
in $B$ by Noetherian induction on $M_4$.
Since the above element is additive on short exact sequences,
it is enough to find a $G$-equivariant submodule
$M$ of $M_4$ which satisfies \eqref{zero}.
Pick an associated prime $p$ of $M_4$ and consider
an element $m\in M_4$ such that ${\rm Ann}(m)=p$. 
In general, we can not expect $m$ to be an eigenelement,
nor can we expect $p$ to be $G$-invariant. However,
let us look at the vector space $V$ which is the linear span
of $gm,g\in G$. Observe that for each $g\in G$, the annihilator
$gp$ of $gm$ is also an associated prime of $M_4$. Since $G$
may not be connected, it could conceivably permute the 
associated primes. Consider the ideal $I\subset A$ given by
$$
I = \bigcap_{g\in G} gp.
$$
It is clear that $I$ is $G$-invariant and that it
annihilates $V$. Hence the submodule $M$ of $M_4$ generated
by $V$ is an $A/I$-module. 

Since $M_4$ is supported on $\CC^n-Z$, each of its associated
primes contains a prime ideal of $A$ that corresponds to
an irreducible component of $\CC^n-Z$. More specifically,
this is an ideal $J$ which is  generated by $z_i$ with indices $i$,
such that $v_i$ do not lie in a cone $C$ of $\Sigma$.
Since $p\supseteq J$, we get $I\supseteq J$, and $M$
is also an $A/J$-module. We now repeat the argument
of Theorem \ref{surj} to resolve $M$ by direct sums of rank one free $G$-equivariant
$A/J$-modules generated by eigenelements. In view
of the long exact sequences of $\Tor$, it suffices to show that 
\eqref{zero} holds true for such rank one $A/J$-module. 
But this is precisely a relation from the second set of relations on $x_i$, times a monomial in $x_i$ to account for possible 
character of the action of $G$ on the generator.

It remains to show that $\rho_1\circ \rho$ is the identity.
Since $B$ is additively generated by monomials 
$\prod_i x_i^{r_i}$, it is enough to check this for 
a monomial. Observe that the global sections of 
$\otimes_i {\mathcal L}_i^{r_i}$ form a free module $A_\chi$ over
$A$ of rank one. Hence $\Tor^{>0}(A_\chi,\CC)$ are zero and 
$\Tor^0(A_\chi,\CC)=\CC$, with character the $\chi$ that corresponds
to $\prod_i R_i^{r_i}$. This finishes the proof.
\end{proof}

\section{Combinatorial Chern character}\label{sec5}
In this section we study the $K$-theory
of the reduced toric DM stack $\PP_{\bf\Sigma}$
in more detail. More specifically, we show that
$K_0(\PP_{\bf\Sigma},\CC):=
K_0(\PP_{\bf\Sigma})\otimes_\ZZ\CC$ 
is isomorphic \emph{as a vector
space} to $H_{SR}(\PP_{\bf\Sigma},\CC)$.

By Theorem \ref{main}, we have that $K_0(\PP_{\bf\Sigma},\CC)\cong
B_\CC:=B\otimes_\ZZ\CC$.
Its maximum ideals correspond to points 
$(y_1,\ldots,y_n)\in \CC^n$ 
that satisfy
$\prod_{i=1}^n y_i^{f(v_i)}=1$ and $\prod_{i\in I}(1-y_i)=0$,
for $f$ and $I$ in the definition of $B$.
\begin{lemma}\label{maxid}
The 
ring $B_\CC$ is Artinian. Its maximum ideals
are in one-to-one correspondence with the 
elements of ${\rm Box}({\bf\Sigma})$ as follows.
A point $v = \sum_{v_i\in C} \alpha_i v_i$ corresponds to
the $n$-tuple
$(y_1,\ldots,y_n)\in\CC^n$ with $y_i=\ee^{2\pi\ii\alpha_i}$
for $v_i\in C$ and $y_i=1$ otherwise. 
\end{lemma}

\begin{proof}
We need to solve 
for $(y_1,\ldots,y_n)\in \CC^n$ that satisfy
$\prod_{i=1}^n y_i^{f(v_i)}=1$ and $\prod_{i\in I}(1-y_i)=0$
as in the definition of $B$. Because of the second set of
equations, $y_i$ are equal to $1$ for all $v_i$ outside some
cone $C\in\Sigma$. We can assume that $C$ is generated
precisely by $v_i$ for indices $i$ with $y_i\neq 1$. 
To simplify notations, let us assume that these $v_i$ 
are $v_1,\ldots,v_k$ for some $k\leq \rk N$.

The first set of relations now reads 
\begin{equation}\label{yrel}
\prod_{i=1}^k y_i^{f(v_i)}=1
\end{equation}
for any linear function $f:N\to \ZZ$. Let $N_1$ be the
intersection of $N$ and the \emph{rational} 
span of $v_1,\ldots,v_k$.
 It is enough to look at $f:N_1\to \ZZ$.
By looking at some $f_i$ which is zero for $j\in[1,\ldots,k]-\{i\}$,
we conclude that $y_i$ is a root of $1$. We introduce
$\alpha_i\in[0,1)$, such that $y_i=\ee^{2\pi\ii \alpha_i}$.
Then the relations \eqref{yrel} amount to
$\sum_i f(v_i)\alpha_i\in\ZZ$ for all
$f\in N_1^\dual$. This is true if and only if
$v=\sum_{i=1}^k\alpha_i v_i\in N_1$.
Hence, the solutions to \eqref{yrel} are in 
one-to-one correspondence with elements of ${\rm Box}(C)$.
The condition $y_i\neq 1$ for $v_i\in C$ assures that
$v$ does not lie in ${\rm Box}(C_1)$ for any proper
face $C_1$ of $C$.  

Since we are looking at all possible cones $C$ here,
the description of maximum ideals follows. Finally,
the ring $B_\CC$ is Artinian, since it is Noetherian
of Krull dimension zero.
\end{proof}

Since $B_\CC$ is Artinian, it is 
a direct sum of Artinian local rings obtained by localizing 
at maximum ideals, which we denote by $(B_\CC)_v$.
We have 
$$
B_\CC=\oplus_{v\in{\rm Box}({\bf\Sigma})}
(B_\CC)_v
$$
The next lemma describes the structure of $(B_\CC)_v$.

\begin{lemma}\label{Ksector}
Let $C$ be the minimum cone of $\Sigma$ that contains 
$v$.  Then the ring $(B_\CC)_v$ is isomorphic as a 
$\CC$-algebra to the
SR-cohomology with complex coefficients
of the (closed) subvariety of $\PP_\Sigma$ that corresponds 
to $C$.
\end{lemma}

\begin{proof}
To simplify notations, we assume that 
$v=\sum_{i=1}^k \alpha_i v_i$ with $\alpha_i\in (0,1)$.
We will also index the rest of $v_i$ in such a way that 
$v_{k+1},\ldots,v_l$ are contained in some cone $C_1\supset C$,
and $v_{l+1},\ldots, v_n$ are not. 

We can localize first and then apply our relations. In fact,
since $B_\CC$ is Artinian,
we may assume to be working in the quotient of 
the power series ring in $x_i-y_i$ by a sufficiently high
power of the maximum ideal. This makes 
$x_i-1$ nilpotent in $(B_\CC)_v$
for $i>k$ and it makes $x_i-\ee^{2\pi\ii\alpha_i}$  nilpotent
for $1\leq i\leq k$. We 
define $z_i=\log(x_i):=\sum_{m>0}\frac 1m(x_i-1)^m$
for $i>k$ and $z_i=\log(x_i\ee^{-2\pi\ii\alpha_i}):=
\sum_{m>0}\frac 1m(x_i\ee^{-2\pi\ii\alpha_i}-1)^m$
for $i=1,\ldots,k$. The elements $z_i$ are also nilpotent and
we can assume to be working in the quotient $B_1$ of 
$\CC[[z_1,\ldots,z_n]]$ by a sufficiently high power of 
the maximum ideal.

We further observe that $z_j=0$ in $(B_\CC)_v$ for $j>l$.
Indeed, we have 
$$
(x_j-1)\prod_{i=1}^k (x_i-1) = 0
$$
which translates into 
$$
(\ee^{z_j}-1)\prod_{i=1}^k (\ee^{2\pi\ii\alpha_i }\ee^{z_i}-1)=0.
$$
Since $\alpha_i\in (0,1)$, we see that $(\ee^{2\pi\ii\alpha_i}\ee^{z_i}-1)$ is invertible in $B_1$,
so $\ee^{z_j}-1=0$. This gives
$z_j(1+\frac 12 z_j+...)=0$, which leads to $z_j=0$. 
As a result, we may just ignore
$z_j$ for $j>l$ in our calculations and work in the 
quotient $B_2$ of 
$\CC[[z_1,\ldots,z_l]]$ by a sufficiently high power of 
the maximum ideal.

The relations $\prod_{i\in I}(x_i-1)=0$ are now 
only nontrivial for $I\subseteq[1,\ldots,l]$. Then every
such set can be enlarged by adding all of $[1,\ldots,i]$.
Consider the quotient fan $\Sigma_C$ in $N/N_1$ which is
made from the images of the cones that contain $C$. 
Since $(x_i-1)$ is invertible for $i\leq k$,
the relations $\prod_{i\in I}(x_i-1)=0$ become relations 
of the form $\prod_{i\in I_1}(\ee^{z_i}-1)$ for $I_1\subseteq
[k+1,\ldots,l]$ such that $z_i,i\in I_1$ are not contained
in any cone of $\Sigma_C$. As before, this is equivalent
to $\prod_{i\in I_1}z_i=0$ for these $I_1\subseteq
[k+1,\ldots,l]$. This completely describes the second set
of relations on $z_i$.

We now need to rewrite the first set of relations in terms of $z_i$.
Let $f:N\to \ZZ$ be a linear function. Then we have
$$
\prod_{i=1}^{k} (\ee^{2\pi\ii\alpha_i})^{f(v_i)}
\prod_{i=1}^{k+l} \ee^{z_if(v_i)}-1=0.
$$
Since $f(v)\in\ZZ$, this becomes simply 
$$
\ee^{\sum_{i=1}^{k+l} z_if(v_i)}-1=0,
$$
which is further equivalent to  
$$
\sum_{i=1}^{k+l} z_i f(v_i) = 0.
$$

We can now ignore the integrality condition on $f$ and 
simply look at all linear functions $f:N\to \ZZ$.
Consider the subspace $V$ of linear functions on $N$ that 
satisfy $f(v_i)=0$ for all $i\leq k$, and let $V_1$ be 
a complement of it. We can use elements of $V_1$ to
express $z_1,\ldots,z_k$ as linear combinations of 
$z_{k+1},\ldots,z_l$. This will allow us to write 
$(B_\CC)_v$ as the quotient of the ring 
$B_3=\CC[z_{k+1},\ldots,z_l]$ by a high power of a maximum
ideal, and by the relations
\begin{itemize}
\item
$\prod_{i\in I_1} z_i$, for $I_1$ such that $z_i,i\in I_1$
do not lie in a cone of $\Sigma_C$,
\item
$\sum_{i=1}^{k+l} z_i f(v_i)$ for any linear function
$f:N/N_1\to \QQ$.
\end{itemize}
This is immediately recognized as the SR-cohomology ring
of the toric subvariety of $\PP_\Sigma$ that corresponds
to the cone $C$.
\end{proof}

\begin{theorem}\label{cherncharacter}
There is a natural vector space isomorphism between
$K_0(\PP_{\bf\Sigma},\CC)$ and  $H_{SR}(\PP_{\bf\Sigma},\CC)$.
\end{theorem}

\begin{proof}
Follows directly from Theorem \ref{main},
Lemmas \ref{maxid} and \ref{Ksector} and 
Propositions \ref{split} and \ref{sector}.
\end{proof}

\begin{remark}
We call the map of Theorem \ref{cherncharacter} \emph{combinatorial Chern character}.
While it is an isomorphism of vector spaces,
we stress that this map is not a ring homomorphism, 
since $B_\CC$ is semilocal and $H_{SR}(\PP_{\bf\Sigma},\CC)$
is local. One motivation behind our construction
is that it generalizes the Chern character for projective
toric varieties (which, however, is a ring isomorphism).
\end{remark}

\begin{remark}
If $\PP_\Sigma$ is  projective, one can use the isomorphism
of SR- and orbifold cohomology to construct
a map $K_0({\bf\Sigma},\CC)\to H_{orb}(\PP_{\bf\Sigma},\CC)$.
We suspect that this map
is  the Chern character map 
of \cite{Adem-Ruan}, which also motivated our terminology.
But since our technique is quite
different, we found it hard to make the connection explicit.
\end{remark}

\begin{remark}
It is an interesting question as to under what conditions
on the fan the $K$-theory is torsion-free. We do not know 
the answer to this, even in the case of smooth toric 
\emph{varieties}.
\end{remark}

\begin{remark}
The additive isomorphism between $K_0(\PP_\Sigma,\CC)$ 
and the SR-cohomology 
$H_{SR}(N,\Sigma,\CC)$ indicates
that the $K$-theory of a reduced toric DM stack possesses
an alternative product structure that makes it into a local
ring. One wonders what the geometric meaning
of this structure may be. A related open question is whether
there is a structure like this for the $K$-theory of
any DM stack, not necessarily
a toric one.
\end{remark}

\section{The $K$-theory of nonreduced toric DM
stacks}\label{sec6}
In this section we extend the calculation of $K$-theory
of toric stacks to the nonreduced case. Since this is 
not the main focus of this paper, we only sketch
the changes necessary to extend the results.

The principal feature of the nonreduced case is that $N$
is now just a finitely generated abelian group, i.e. it is allowed
to have torsion. A fan in $N$ is a pullback of a fan in $N/torsion$.
A stacky fan is defined by a choice of a nonzero element $v_i$ in
each of the one-dimensional cones of $\Sigma$. We still have
a map 
$$
\ZZ^n\to N
$$
with finite coindex. This allows one to define a Gale dual
$$
(\ZZ^n)^\dual\to N'
$$
which is the analog of $(\ZZ^n)^\dual\to\Coker\phi$ of Section
\ref{sec2},
but we no longer have the surjectivity. 
We refer the reader to \cite{BCS} for the details.
Consequently, the group $G={\rm Hom}(N',\CC)$
is no longer a subgroup of $(\CC^*)^n$, but
rather maps to it with a finite kernel.

Otherwise, the definition of the open subset 
$Z\subseteq \CC^n$ is the same as in the reduced case.
This allows one to define line bundles ${\mathcal L}_i$
on $\PP_{\bf\Sigma}$. Unfortunately, they will no longer generate
the $K$-theory. The problem is that in the proof of Theorem \ref{surj} we have used that every character of $G$ lifts to
a character of $(\CC^*)^n$. This is no longer the case.
However, we can still look at the ring $B$ which is the 
quotient of the character ring of $G$ by the relations
$$
\prod_{i\in I} (x_i-1) = 0 
$$
for all $I\subseteq [1,\ldots,n]$, such that $v_i,i\in I$ are 
not contained in any cone of $\Sigma$. Here $x_i$ correspond
to ${\mathcal L}_i$ as in nonreduced case. With this modification,
the proofs of 
both Theorems \ref{surj} and \ref{main} are extended to
the nonreduced case without any major changes to show 
that $K_0(\PP_{\bf\Sigma})\cong B$.

The description of the maximum ideals of $B_\CC$ from
Lemma \ref{maxid} still holds in the nonreduced case,
but the proof is a bit more complicated, especially
is one does not assume the condition \eqref{technical}.
 Specifically, these
ideals correspond to elements of $G$ whose action
on $\CC^n$ has eigenvalues one outside of the set of
indices $i$ for $v_i$ in some cone $C\in \Sigma$.
So we have a morphism $\psi:N'\to \CC^*$ such that the 
composition $(\ZZ^n)^\dual\to N'\to \CC^*$ takes value
one on the basis elements that correspond to $v_i$
outside of $C$. The value of $\psi$ on any element of $N'$
is a root of one. Indeed, $\psi\in G$ fixes a point in $Z$,
and it has been shown in \cite{BCS} that all isotropy subgroups
of $G$-action are finite. Consequently, we can think of $\psi$
as a map $N'\to \QQ/\ZZ$. Conversely, every such $\psi$
gives rise to a group element with the above eigenvalue
properties, and hence to a local subring of $B_\CC$.

We split
$(\ZZ^n)^\dual$ into $(\ZZ^{\dim C})^\dual\oplus
(\ZZ^{n-\dim C})^\dual$ according to whether the corresponding
basis elements lie in $C$. Consider the short exact
sequence of complexes 
$$
\begin{array}{ccccccc}
0\to&(\ZZ^{n-\dim C})^\dual&\to&(\ZZ^{\dim C})^\dual\oplus
(\ZZ^{n-\dim C})^\dual&\to&(\ZZ^{\dim C})^\dual&\to0\\
&\downarrow&&\downarrow&&\downarrow&\\
0\to&N'&\to&N'&\to&0&\to0
\end{array}
$$
where we assume that the vertical lines are extended by zeroes
in both directions. Gale duality
with torsion (see \cite{BCS}) implies that $N$ can be canonically
identified with the $H^1$ of the derived $Hom(*,\ZZ)$ of the middle 
complex. The above exact sequence shows that 
$N/Span(v_i\in C)$ is $H^1$ of the derived 
$Hom((\ZZ^{n-\dim C})^\dual\to N',\ZZ)$. 

We observe that $\psi: N'\to \QQ/\ZZ$ constructed earlier
can be identified with elements of $H^0$ of 
derived $Hom((\ZZ^{n-\dim C})^\dual\to N',\QQ/\ZZ)$.
The short exact sequence $0\to \ZZ\to \QQ\to\QQ/\ZZ\to 0$
yields an exact sequence
$$
H^0(Hom((\ZZ^{n-\dim C})^\dual\to N',\ZZ))\to H^1(
Hom((\ZZ^{n-\dim C})^\dual\to N',\ZZ))\to
$$
$$
\to H^1(
Hom((\ZZ^{n-\dim C})^\dual\to N',\QQ)).
$$
Since the last term of the above sequence  is torsion free, the group $H^0(Hom((\ZZ^{n-\dim C})^\dual\to N',\ZZ))$
is precisely the torsion subgroup of $N/Span(v_i\in C)$.
Finally, elements of this torsion subgroup are in one-to-one
correspondence with elements of ${\rm Box}(C)$. We
leave to the reader to check that these identifications
are compatible with embeddings ${\rm Box}(C_1)\subseteq 
{\rm Box}(C_2)$ for $C_1\subseteq C_2$ and that they
coincide with the construction of Section \ref{sec5} in
the reduced case.

The 
SR-cohomology of $\PP_{\bf\Sigma}$ is again
given as the quotient of $\CC[N,\Sigma]$ by the linear
relations $\sum_i f(v_i)[v_i]$. Note that it is no longer
a local ring, because of the graded zero part
is the group ring of the torsion subgroup of $N$. 
SR-cohomology again splits into the untwisted and 
twisted sectors according to elements
of ${\rm Box}(\bf\Sigma)$, and Propositions \ref{split}
and \ref{sector} still hold.

The combinatorial
Chern character map of Theorem \ref{cherncharacter}
 still exists in
the nonreduced case. It sends the local subring of the 
Artinian ring $B_\CC$ which corresponds to the 
point $v\in{\rm Box}({\bf\Sigma})$ to the 
SR-cohomology of the 
corresponding twisted sector. More specifically, let
$v$ be this point and let $C$ be the minimum cone of $\Sigma$
that contains it. The map from $G$ to 
$G' $ is a finite unramified cover, so the local ring summand for
$B_\CC$ is isomorphic to the local ring of a semi-local
ring $B'_\CC$
which is obtained by cutting the character ring of $G'$ by
the same set of relations. 
Then the calculation of Lemma \ref{Ksector} shows
that $(B_\CC)_v\cong (B'_\CC)_v$ is isomorphic to
the SR-cohomology of the twisted sector of the reduction
of $\PP_{\bf\Sigma}$ that corresponds
to the reduction of $v$ modulo torsion. However,
this reduction does not change the coarse moduli space of 
the twisted sector, which finishes the construction of 
the isomorphism. The details are left to the reader.

\section{Birational morphisms of reduced toric DM stacks}\label{sec7}
The statements of this section are implicit in \cite{Kawamata},
but we adjust the exposition for our notations.

Let $N$ be a lattice and let ${\bf\Sigma}=(\Sigma,\{v_i\})$ 
be a stacky fan
in $N$. We will assume that every cone of $\Sigma$ is contained in a cone of $\Sigma$ of dimension ${\rm rk} N$. 
As before, let $n$ be the number of $v_i$ and let
$Z\subset\CC^n$ consist of the
points ${\bf z}=(z_1,\ldots,z_n)$ such that the set of $v_i$ for
the zero coordinates of $\bf z$ is contained in a cone of $\Sigma$.
Let $G$ be the subgroup of $(\CC^*)^n$ described by 
\eqref{eqlam}.
The toric DM stack $\PP_{\bf\Sigma}$
is defined as the stack quotient 
$[Z/G]$ where $Z$ and $G$ are 
endowed with the natural reduced 
scheme structures and the action of $(\lambda_i)$ on
$(z_i)$ is given by $(\lambda_iz_i)$.

Let ${\bf\Sigma'}=(\Sigma',\{v'_j\})$ be another stacky fan in the same lattice
$N$. Let us further assume that every cone of $\Sigma'$
is contained in a cone of $\Sigma$. 
In particular, for every $j$ there is a unique way of writing
$v'_j$ in the form
$$
v'_j = \sum_i \alpha_{i,j} v_i, ~\alpha_{i,j}\in \QQ
$$
under the assumption that $\alpha_{i,j}$ are zero unless 
$v_i$ lies in the minimum cone of $\Sigma$ that contains 
$v'_j$. Moreover, let us assume
that all $\alpha_{i,j}$ above are \emph{integer}.
\begin{definition}\label{birmor}
Under the above assumptions, there is a morphism
$$
\PP_{\bf\Sigma'}\to \PP_{\bf\Sigma }
$$ 
defined as follows. Define a homomorphism 
$(\CC^*)^{n'}\to (\CC^*)^n$
which sends $(\lambda'_j)$ to $(\prod_j (\lambda'_j)^{\alpha_
{i,j}})$. It is easy to see that this homomorphism sends points
of $G'$ to points of $G$. We also consider the map from
$\CC^{n'}$ to $\CC^n$ which sends $(z'_j)$
to $(\prod_j (z'_j)^{\alpha_{i,j}})$. This map sends points in
$Z'$ to points in $Z$. Indeed, if the $C'$ is a cone of 
$\Sigma'$ that contains all $v'_j$ for  which $z'_j=0$,
then $\prod_j (z'_j)^{\alpha_{i,j}}$ is nonzero unless
$v_i$ lies in the minimum cone of $\Sigma$ that contains
$C'$.
The map $Z'\to Z$ is compatible with the map of groups
$G'\to G$ and their actions. This gives a map of quotient
stacks, since one has  the map between the corresponding 
groupoids.
\end{definition}

\begin{remark}
We call the toric morphisms of Definition \ref{birmor}
\emph{birational}, because they induce isomorphisms
on the big strata $(\CC^*)^{{\rm rk}N}$ in both
stacks. 
\end{remark}

\section{$K$-theory pullbacks for birational morphisms}
\label{sec8}
Let $\mu:\PP_{\Sigma',\{v_j'\}}\to \PP_{\Sigma,\{v_i\}}$ 
be a toric birational morphism of toric DM stacks.
Let $R_i$ be defined as the elements of 
$K_0(\PP_{\Sigma,\{v_i\}},\QQ)$
that correspond to the invertible sheaves that correspond to $v_i$,
and similarly for $R'_j$. There is a pullback map 
$$\mu^*:K_0(\PP_{\Sigma,\{v_i\}},\QQ)\to
 K_0(\PP_{\Sigma',\{v'_j\}},\QQ)$$
defined by pulling back the coherent sheaves.
As we saw earlier in Theorem \ref{main},  the elements
$\prod_i R_i^{r_i}$ span $K_0(\PP_{\Sigma,\{v_i\}},\QQ)$.

\begin{proposition}\label{pull}
The pullback map $\mu^*$ is given by
$$\mu^*\prod_i R_i^{r_i} = \prod_j (R'_j)^{\sum_i
\alpha_{i,j}r_i}$$
\end{proposition}

\begin{proof}
The category of coherent sheaves on $\PP_{\Sigma,\{v_i\}}$ 
is equivalent to that of $G$-linearized coherent sheaves on
$Z$, see \cite[Example 7.21]{Vistoli}. 
Under this equivalence, $\prod_i R_i^{r_i}$ corresponds
to the trivial sheaf on $Z$, linearized 
by 
${\mathcal O}_Z\to g^*{\mathcal O}_Z={\mathcal O}_Z$ 
which send
$1$ to $\prod_i (\lambda_i)^{r_i}$.
These isomorphisms pull back to the isomorphisms
of ${\mathcal O}_{Z'}$
that, for given $(\lambda'_1,\ldots,\lambda'_{n'})$, send $1$ to 
$\prod_{i,j}((\lambda'_j)^{\alpha_{i,j}})^{r_i}$
on $Z'$. Since the power of $\lambda'_j$ is
$\sum_i \alpha_{i,j}r_i$, the result follows.
\end{proof}

\begin{remark}
It is an amusing, though unnecessary, combinatorial
exercise to see that the relations among $\prod_i R_i^{r_i}$
get mapped in the ideal of the relations among $\prod_j (R'_j)^{r'_j}$.
\end{remark}

\section{Weighted blowups and pushforward formulas}
\label{sec9}
Let $(\Sigma,\{v_i\})$ be a stacky fan in lattice $N$.
Let $C$ be a cone in $\Sigma$ of dimension $d>1$ and
let $\{v_1,\ldots v_d\}$ be the set of $v_i\in C$.
Let $h_1,\ldots,h_d$ be some positive integers. 
Consider 
$$v'_0:=\sum_{i=1}^d h_i v_i.$$ 
The cone $C$ decomposes
into a union of $d$ cones of dimension $d$ 
which are given by positive linear combinations of 
$v'_0$ and all but one of $v_i,1\leq i\leq d$.

We define a stacky fan $(\Sigma',\{v'_j\})$ as follows.
We add an extra ray which corresponds to $v'_0$
and keep the rest
of $v_i$ unchanged. Thus $n'=n+1$ and $v'_i=v_i$
for $i>0$.
Every 
cone $\hat C$ of $\Sigma$ that does not contain $C$ 
is still a cone of $\Sigma'$. Every cone $\hat C\supseteq C$
is subdivided into $d$ cones $C'_1,\ldots,C'_d$ according
to the above decomposition of $C$. 

Our definitions assure that there is a birational morphism
$$
\mu:\PP_{\Sigma',\{v_j'\}}\to \PP_{\Sigma,\{v_i\}}
$$
which we call the \emph{weighted blowdown} morphism.
We will also denote by $\pi$ the corresponding
morphism $Z'\to Z$.
Our goal is to calculate the  $K$-theory pushforward 
of $\mu$. It is best described in terms of generating
functions.

\begin{theorem}\label{push}
Let $R=R'_0$ be the $K$-theory class
of the invertible sheaf that corresponds 
to the extra ray $v'_0$.
There holds
$$
\mu_*\Big( \frac 1{1-R^{-1}t} \Big)=\frac 1{1-t}
-\frac t{1-t} \prod_i \frac {1-R_i^{-1}}{1-R_i^{-1}t^{h_i}}
$$
which should be interpreted as an identity of formal
power series in $t$ with values in $K$-theory.
\end{theorem}

\begin{proof}
We will again identify the abelian categories of sheaves
of $\PP_{\Sigma',\{v_j'\}}$ and $\PP_{\Sigma,\{v_i\}}$
with the categories of $G'$- and $G$-linearized sheaves
on $Z'$ and $Z$ respectively. The $K$-theory pushforward
of a sheaf $F'$ is defined as the $K$-theory image of the 
alternating sum of the higher direct images of $F'$ under
$\mu$. In order to describe the higher direct images of 
$\mu$, we need to describe the direct image functor for $\mu$
in terms of $\pi:Z'\to Z$.

We first observe that the kernel $H$ of $G'\to G$ 
is isomorphic to $\CC^*$. Indeed, the map comes from
$(\CC^*)^{n+1}\to (\CC^*)^n$ defined by
$$
(\lambda'_0,\lambda'_1,\ldots,\lambda'_n)
\mapsto 
(\lambda'_1(\lambda'_0)^{h_1},\ldots,\lambda'_d
(\lambda'_0)^{h_d},
\lambda'_{d+1},\ldots,\lambda'_n).
$$
The kernel is given by $\lambda'_0=\lambda$,
$\lambda'_i=\lambda^{-h_i}$ for $1\leq i\leq d$,
$\lambda'_{>d}=1$, which clearly lies in $G'$.
Moreover, the map $G'\to G$ is split surjective. 
Indeed, given $(\lambda_i)\in G$, the collection
$(1,\lambda_1,\ldots,\lambda_n)$ will lie in $G'$,
which produces the splitting. We will
only need surjectivity for our arguments.

For every $G'$-linearized 
coherent sheaf ${\mathcal F}'$ on $Z'$, its pushforward
${\mathcal F}=\pi_*{\mathcal F}'$ 
on $Z$ inherits the $G'$-linearization under
the action of $G'$ on $Z$ induced from $G'\to G$.
Its subsheaf ${\mathcal F}^H$ of $H$-invariants is therefore given
the structure of a $G$-linearized sheaf on $Z$.
The $K$-theory pushforward
of a sheaf on $\PP_{\Sigma',\{v_j'\}}$ 
that corresponds to ${\mathcal F}'$ is given by ${\mathcal F}^H$.
Indeed, consider the following commutative diagram
$$
\begin{array}{rcccl}
&Z'&\stackrel \pi\to&Z&\\
&q'\downarrow~&&~\downarrow q&\\
&[Z'/G']&\stackrel \mu\to&[Z/G]&
\end{array}
$$
The sheaf on $[Z'/G']$ that corresponds to ${\mathcal F}'$ is
the subsheaf of $G'$-invariants of $q'_*{\mathcal F}'$. Consequently,
its direct image in $[Z/G]$ is the subsheaf of $G'$-invariants
of $q_*\pi_* {\mathcal F}'=q_*{\mathcal F}$. 
The subgroup $H$ of $G'$ acts trivially
on $Z$, so the $G'$-invariants of $q_*{\mathcal F}$ 
are the $G$-invariants of $q_*{\mathcal F}^H$. This 
corresponds to the $G$-equivariant sheaf ${\mathcal F}^H$ on $Z$.

We have thus described the direct image functor in terms
of the composition of the direct image functor $\pi_*$ from
$Coh_{G'}(Z')$ to $Coh_{G'}(Z)$ and the functor of $H$-invariants
from $Coh_{G'}(Z)$ to $Coh_{G}(Z)$. The latter is exact, therefore,
the higher direct images of $\mu$ are given by 
$(R^k\pi_*({\mathcal F}'))^H$.

In order to prove the theorem, we need to calculate 
the higher direct images of the 
$G'$-linearized invertible sheaf $\mathcal L$ on $Z'$ 
with the linearization 
that sends $1$ to $(\lambda'_0)^{-l}$ for some integer $l\geq 0$.
We claim
that  
$$(R^{>0}\pi_* {\mathcal L})^H=0.$$
It is sufficient 
to check the statement at a fiber. Let 
${\bf z}=(z_1,\ldots, z_n)$ be a point in $Z$. The fiber of $\pi$
consists of points $(z'_0,z'_1,\ldots,z'_n)\in Z'$ such
that 
$$
{\bf z}= (z'_1(z'_0)^{h_1},\ldots,z'_d
(z'_0)^{h_d},
z'_{d+1},\ldots,z'_n).
$$
Hence if one or more of $z_1,\ldots,z_d$ are nonzero,
then the fiber is isomorphic to $\CC^*$, which is an orbit
of $H$. The structure sheaf of $\CC^*$ has no higher 
cohomology, so $R^{>0}\pi_*{\mathcal L}$ 
will have zero fibers at such $\bf z$,
and $R^{>0}\pi_*{\mathcal L}$ are supported over $z_1=\ldots=z_d=0$.
For these $\bf z$, the fiber is described by $z'_0=0$,
$z'_i=z_i$ for $i>d$, and $(z'_1,\ldots,z'_d)\neq(0,\ldots,0)$.
The fiber is isomorphic to $\CC^d-{\bf 0}$, with the group
$H$ acting by multiplications of the $i$-th coordinate 
by $\lambda^{-h_i}$. The cohomology of 
the structure sheaf ${\mathcal O}$ on $\CC^d -{\bf 0}$
occurs at $H^0$ and $H^{d-1}$  only, as can be easily
calculated by the \v{C}ech complex for the covering
by the open sets $(z'_i\neq 0)$. Moreover, there is a natural
isomorphism 
$$H^{d-1}(\CC^d -{\bf 0},{\mathcal O})\cong 
(\prod_i (z'_i)^{-1}) \CC[(z'_1)^{-1},\ldots,(z'_d)^{-1}].$$ 
This means that the action of $\lambda\in H$ multiplies
monomial generators 
 of $H^{d-1}(\CC^d -{\bf 0},{\mathcal O})$ by
$\lambda^{-\sum_{i}h_is_i}$ for some $s_i<0$. Consequently,
it multiplies generators of $H^{d-1}(\CC^d -{\bf 0},{\mathcal L})$
by $\lambda^{l-\sum_i h_i s_i}$, which shows
that the space of $H$-invariants is zero.
Of course, the above is basically a calculation of cohomology
of a line bundle on a weighted projective space.

As a result, to calculate the $K$-theory pushforward of 
$R^{-l}$ for $l\geq 0$ it is enough to calculate the direct image
of the corresponding invertible sheaf.
In the notations of the preceding paragraph, 
$\mathcal L$ is naturally embedded into $\mathcal O$
as an ideal sheaf of $(z'_0)^{l}$. As a result, $(\pi_*{\mathcal L})^H$ is
embedded into $(\pi_*{\mathcal O})^H$. 
We  first show
that the latter is isomorphic to ${\mathcal O}_{Z}$, which
is just the statement that $\mu_*{\mathcal O}={\mathcal O}$
for a birational morphism. 

The isomorphism will be glued from the isomorphisms on
$(\CC^*)^n$-invariant affine subsets $U_\sigma
\subseteq Z$
which correspond to the maximum-dimensional cones $\sigma$
of $\Sigma$.
The subset $U_\sigma$ is defined by the condition that 
$z_i\neq 0$ for all $v_i\not\in \sigma$ and is isomorphic to
$(\CC^*)^{n-{\rm rk}N}\times \CC^{{\rm rk} N}$.
If $\sigma\not\supseteq C$, then the
preimage $\pi^{-1}U_\sigma$ in $Z'$ is given by the conditions
$z'_0\neq 0$ and $z'_i\neq 0$ for $i>0$ and $v_i\not\in\sigma$.
Then $\pi^{-1}U_\sigma$ is isomorphic to $(\CC^*)^{n+1-{\rm rk}N}\times \CC^{{\rm rk} N}$. The sections of 
$\pi_*{\mathcal O}_{Z'}$
on $U_\sigma$ are spanned by the monomials $(z'_0)^{s_0}\prod_i (z'_i)^{s_i}$ with $s_i\geq 0$ for all 
$v_i\in\sigma$. The $H$-invariant sections in addition satisfy
$s_0=\sum_{i=1}^d h_i s_i$. The isomorphism from
${\mathcal O}_{U_\sigma}$ to $(\pi_*{\mathcal O}_{Z'})^H$
is constructed by sending $\prod_{i=1}^n z_i^{s_i}$ to
$(z'_0)^{\sum h_i s_i}\prod_{i=1}^n (z'_i)^{s_i}$. In the case
when $\sigma\supseteq C$, the preimage 
$\pi^{-1}U_\sigma$ is given by the conditions
$z'_i\neq 0$ for $v_i\not\in\sigma$, and $(z'_1,\ldots, z'_d)
\neq {\bf 0}$. Consequently, it is isomorphic to
$\CC\times(\CC^d-{\bf 0})
\times (\CC^*)^{n-{\rm rk}N}\times\CC^{{\rm rk}N-d}$.
Because $d\geq 2$, the sections of $\mathcal O$
on this space ignore the deletion of $\bf 0$.
The sections of $\pi_*{\mathcal O}_{Z'}$ are over
$U_\sigma$ are therefore spanned
by monomials $(z'_0)^{s_0}\prod_i (z'_i)^{s_i}$ 
with $s_0\geq 0$ and $s_i\geq 0$ for $v_i\in\sigma$.
The $H$-invariance condition simply expresses $s_0$
in terms of other $s_i$ but imposes no further restrictions
on $s_1,\ldots,s_n$. Consequently, we again have 
an isomorphism between ${\mathcal O}_Z$
and $(\pi_*{\mathcal O}_{Z'})^H$. It is clear that these
isomorphisms are compatible on the intersections
and hence glue together to show $(\pi_*{\mathcal O}_{Z'})^H
\cong {\mathcal O}_Z$.

Obviously, the sheaf ${\mathcal O}_Z$
is a restriction of ${\mathcal O}_{\CC^n}$ via the open
embedding $Z\subseteq\CC^n$. We claim that 
$(\pi_*{\mathcal L})^H$ is also a restriction of an ideal sheaf 
$\mathcal I$ from $\CC^n$. Namely, consider the ideal
sheaf  $\mathcal I$ on $\CC^N$ which corresponds
to the submodule over $\CC[z_1,\ldots,z_n]$ which is the
span of monomials $\prod_{i=1}^n z_i^{s_i}$ with
$$
\sum_{i=1}^d s_i h_i \geq l.
$$
To show that it restricts to $(\pi_*{\mathcal L})^H$ on $Z$,
it is again enough to calculate the sections over the 
open subsets $U_\sigma$. For $\sigma\not\supseteq C$,
we have ${\mathcal O}_Z\vert_{U_\sigma}\cong
{\mathcal I}\vert_{U_\sigma}$. For $\sigma\supseteq C$,
the condition on $(z'_0)^{s_0}\prod_i (z'_i)^{s_i}$
to lie in the ideal of $(z'_0)^l$ translates into $s_0\geq l$.
The invariants of that come from the monomials
$\prod_i z_i^{s_i}$ with $\sum_{i=1}^d s_i h_i \geq l$,
as claimed.

To calculate the pushforward of $\mathcal L$ 
we now simply need to calculate the free graded resolution of 
the ideal $I$ of the polynomial ring $\CC[z_1,\ldots,z_n]$
which is the span of the monomials with the condition
$\prod_i z_i^{s_i}$ with $\sum_{i=1}^d s_i h_i \geq l$.
Clearly, the variables $z_{d+1},\ldots,z_n$ can be ignored.
If 
$$
0\to F^d\to F^{d-1}\to\cdots\to F^0\to I\to 0
$$
is such a free resolution, then the alternating sum of the 
$\ZZ^d$-graded dimensions of $F^k$ is the 
$\ZZ^d$-graded dimension of $I$, and the same is true
for their generating functions. The generating function
of a copy $A$ of $\CC[z_1,\ldots,z_n]$ with the grading shifted
so that the multidegree of $1$ is $(r_1,\ldots,r_d)$ 
is 
$$
\sum_{{\bf deg}\in\ZZ^d}\dim_\CC(A_{\bf deg}){\bf t}^{\bf deg}=
\prod_{i=1}^d \frac {t_i^{r_i}}{1-t_i}.
$$
On the the other hand, such a module $A$ gives rise to the 
invertible sheaf on $Z$ that gives $\prod_{i=1}^d R_i^{-r_i}$
in $K$-theory
of $\PP_{\Sigma,\{v_i\}}$.
So the $K$-theory pushforward of $R^{-l}$ is given by
$$
\mu_*R^{-l}=\prod_{i=1}^d(1-R_i^{-1})
\sum_{(s_1,\ldots,s_d)\in \ZZ_{\geq 0}^d,~
h_1s_1+\cdots+h_ds_d\geq l
} R_1^{-s_1}\cdot\ldots\cdot R_d^{-s_d}
$$
where the sum should be interpreted as a formal power series
in $R_i^{-1}$ which actually gives a polynomial after being
multiplied by $\prod (1-R_i^{-1})$.

It remains to observe that if we look at the above in terms of the
generating functions for all $l\geq 0$, then 
$$
\sum_{l\geq 0} t^l\mu_*R^{-l} 
=
\prod_{i=1}^d(1-R_i^{-1})
\sum_{l\geq 0} t^l\sum_{(s_1,\ldots,s_d)\in \ZZ_{\geq 0}^d,~
\sum_{i=1}^d h_is_i\geq l
} \prod_{i=1}^dR_i^{-s_i}
$$
$$
=\prod_{i=1}^d(1-R_i^{-1})
\sum_{(s_1,\ldots,s_d)\in \ZZ_{\geq 0}^d}
\prod_{i=1}^dR_i^{-s_i}
\sum_{0\leq l\leq \sum_{i=1}^d h_is_i} t^l
$$
$$
=\prod_{i=1}^d(1-R_i^{-1})
\sum_{(s_1,\ldots,s_d)\in \ZZ_{\geq 0}^d}
\frac {1-t^{1+\sum_i h_is_i}}{1-t}
\prod_{i=1}^dR_i^{-s_i} 
$$
$$
=\frac 1{1-t}-
\frac t{1-t}\prod_{i=1}^d(1-R_i^{-1})
\sum_{(s_1,\ldots,s_d)\in \ZZ_{\geq 0}^d}
{t^{\sum_i h_is_i}}\prod_{i=1}^dR_i^{-s_i}
$$
$$
=\frac 1{1-t}-
\frac t{1-t}\prod_{i=1}^d\frac{1-R_i^{-1}}
{1-R_i^{-1}t^{h_i}}.
$$
We remark that the above calculations should be interpreted
as calculations in formal power series in $t$ and $R_i^{-1}$
with only finitely many terms at any given degree, so 
convergence is never an issue.
This gives the desired formula for 
$\mu_*\frac 1{1-R^{-1}t}$.
\end{proof}

\begin{remark}
The above theorem allows one to calculate the pushforward
of any element of $K$-theory in view of the description of 
the pullback in Proposition \ref{pull}
and the pull-push formula. Indeed, every 
element of $K_0(\PP_{\Sigma',\{v'_i\}})$ can be written 
as a polynomial in $R$ and $R^{-1}$ with coefficients 
in $\mu^*K_0(\PP_{\Sigma,\{v_i\}})$. Since $R$ is
quasi-unipotent, it can be expressed in terms 
of negative powers of $R$.
\end{remark}

So far we have considered the weighted blowups
of a stratum of codimension $d>1$. While blowups in codimension
one are isomorphisms in the smooth variety case, this is no
longer true for weighted blowups of stacks.
Namely, let $(\Sigma,\{v_i\})$ be a stacky fan in $N$. Let
$C$ be a dimension one cone of $\Sigma$ and let $v_1$
be the chosen lattice point on $C$. For a positive integer $k$
consider the stacky fan $(\Sigma',\{v'_i\})$ defined by 
$\Sigma'=\Sigma$, $v'_1=kv_1$, $v'_i=v_i, i>1$.
There is a birational  morphism
$$
\mu:\PP_{\Sigma',\{v'_i\}}\to \PP_{\Sigma,\{v_i\}}
$$
which we again call weighted blowdown morphism.
To complete the discussion, we calculate the corresponding pushforward in $K$-theory. Since the pullback is given
by $\mu^*R_1=R_1^k$, $\mu^*R_i=R'_i$ for $i>1$, it is
enough to calculate the pushforward of $R_1^{m}$ for 
$m=1,\ldots, k$. 
\begin{proposition}\label{pushzero}
For $1\leq m\leq k$ we have
$\mu_*(R'_1)^{-m} = R_1^{-1}$. 
\end{proposition}

\begin{proof}
The map $G'\to G$ is given by 
$$(\lambda'_1,\lambda'_2,\ldots,\lambda'_n)
\to(\lambda'_k,\lambda'_2,\ldots,\lambda'_n)
$$
and the map $\pi:Z'\to Z$ is given by 
$$
\pi:(z_1,z_2,\ldots,z_n)\mapsto (z_1^k,z_2,\ldots,z_n).
$$
We denote by $H$ the kernel of $G'\to G$.
Since $\pi$ is finite, its higher direct images vanish, and 
we only need to calculate
$(\pi_* {\mathcal L})^H$ for the sheaf $\mathcal L$ which is
isomorphic to ${\mathcal O}_{Z'}$ with linearization given by
multiplication by $(\lambda'_1)^{-m}$.
We can think of $\mathcal L$ as an ideal sheaf of 
${\mathcal O}_{Z'}$ generated by $(z'_1)^m$. 

Similar to the proof of Theorem \ref{main}, 
$(\pi_* {\mathcal L})^H$  is an ideal sheaf in ${\mathcal O}_{Z}$. 
It is induced from a sheaf on $\CC^n\supset Z$. A monomial
$\prod_i z_i^{s_i}$ lies in the corresponding ideal of 
$\CC[z_1,\ldots,z_n]$ if and only if 
$s_1k\geq m$. Since $1\leq m\leq k$,
this is equivalent to $s_1\geq 1$. This ideal sheaf 
is then identified with ${\mathcal L}_1^{-1}$.
\end{proof}

We rewrite the result of Proposition \ref{pushzero}
to resemble that of Theorem \ref{push}.
We denote $R=R_1'$.
\begin{corollary}
For the blowdown of codimension one,
$$
\mu_*\Big( \frac 1{1-R^{-1}t} \Big)=
\frac 1{1-t} - \frac t{(1-t)} \frac {
(1-R_1^{-1})}
{(1-R_1^{-1}t^k)}.
$$
\end{corollary}

\begin{proof}
$$
\mu_*\Big( \frac 1{1-R^{-1}t} \Big)=
1+\sum_{l>0} \mu_*R^{-l}t^l
=1 + \sum_{l\geq 0} \sum_{m=1}^{k}\mu_*R^{-lk-m}t^{lk+m}
$$
$$
=1 + \sum_{l\geq 0} \sum_{m=1}^{k}\mu_*
(R^{-m}\mu^*R_1^{-l})t^{lk+m}=
1 + \sum_{l\geq 0} \sum_{m=1}^{k}R_1^{-l-1}t^{lk+m}
$$
$$
= 1 + 
\frac {R_1^{-1}(t-t^{k+1})}
{(1-R_1^{-1}t^k)(1-t)} 
=\frac 1{1-t} - \frac t{(1-t)} \frac {
(1-R_1^{-1})}
{(1-R_1^{-1}t^k)}.
$$
\end{proof}

\end{document}